\documentclass[11pt]{amsart}

\usepackage{amsmath,amsthm,amssymb,amsfonts,ifpdf}
\usepackage{xcolor}

\usepackage{xargs}  
\usepackage[colorinlistoftodos,prependcaption,textsize=tiny,obeyFinal]{todonotes}
\newcommandx{\ebltodo}[2][1=]{\todo[linecolor=red,backgroundcolor=red!25,bordercolor=red,#1]{#2}}
{
  \color{olive}%
}%
{}
\PassOptionsToPackage{hyphens}{url}
\ifpdf
\usepackage[pdftex,pdfstartview=FitH,pdfborderstyle={/S/B/W 1},%
colorlinks=true, linkcolor=blue, urlcolor=blue, citecolor=blue,%
pagebackref=true]{hyperref}
\else
\usepackage[hypertex]{hyperref}
\fi
\usepackage{enumitem}
\setenumerate{leftmargin=*}
\setitemize{leftmargin=*}
\usepackage{amscd}
\usepackage[latin1]{inputenc}
\DeclareMathAlphabet{\mathpzc}{OT1}{pzc}{m}{it}
\usepackage{bbm}

\numberwithin{equation}{section}

\swapnumbers
\newtheorem{thm}[subsection]{Theorem}
\newtheorem{coro}[subsection]{Corollary}
\newtheorem*{cor*}{Corollary}

\newtheorem*{thm*}{Theorem}
\newtheorem*{thma*}{Theorem A}
\newtheorem*{thmb*}{Theorem B}
\newtheorem*{thmc*}{Theorem C}

\swapnumbers
\theoremstyle{definition}

\newcounter{consta}
\renewcommand{\theconsta}{{A_{\arabic{consta}}}}
\newcommand{\consta}{\refstepcounter{consta}\theconsta}

\newcounter{constk}
\renewcommand{\theconstk}{{\kappa_{\arabic{constk}}}}
\newcommand{\constk}{\refstepcounter{constk}\theconstk}

\newcounter{constc}

\newcounter{constE}

\newcounter{constd}

\makeatletter
\newcommand*\bigcdot{\mathpalette\bigcdot@{.5}}
\newcommand*\bigcdot@[2]{\mathbin{\vcenter{\hbox{\scalebox{#2}{$\m@th#1\bullet$}}}}}
\makeatother

\def\XXint#1#2#3{{\setbox0=\hbox{$#1{#2#3}{\int}$ }
\vcenter{\hbox{$#2#3$ }}\kern-.6\wd0}}

\DeclareMathOperator{\diff}{d}

\DeclareMathOperator\Ad{Ad}

\newcommand\vol{{\rm{vol}}}
\newcommand\SL{{\rm{SL}}}
\newcommand\PSL{{\rm{PSL}}}

\newcommand\SO{{\rm{SO}}}
\newcommand\Lie{{\rm Lie}}

\def\bbz{\mathbb{Z}}

\def\bbr{\mathbb{R}}

\def\Z{\bbz}
\def\R{\bbr}

\def\hfrak{\mathfrak{h}}

\def\rfrak{\mathfrak{r}}

\def\gfrak{\mathfrak{g}}

\def\vare{\varepsilon}

\def\zg0{Z_{G_\omega}(s)}

\def\zg{Z_G(s)}

\def\be{\begin{equation}}
\def\ee{\end{equation}}

\def\dist{{\rm dist}}

\def\Sob{{\mathcal S}}

\def\dist{d}

\def\rwm{\nu}

\def\rel{r}

\def\ave{\int_{0}^1}
\def\uvk{u_\rel}
\def\uvkd{\diff\!\rel}

\def\mfm{f}

\def\margI{I}

\def\inj{{\rm inj}}

\def\nuni{e}
\def\coneH{\mathsf E}
\def\cone{\mathcal E}

\newcommand{\rhsc}{b}

\newcommand{\scmf}{b}

\newcommand{\trct}{\mathsf M}

\newcommand\egbd{\Upsilon}
\newcommand\pvare{\mathsf{c}}

\newcommand {\absolute}[1] {\left| {#1} \right|}
\newcommand {\norm}[1] {\left\| {#1} \right\|}

\newcommand{\hide}[1]{}

\newcommand\sqf{Q_0}
\newcommand{\cox}{{\mathsf x}}

\title[Polynomially effective Oppneheim]{An effective version of the Oppenheim conjecture with a polynomial error rate}

\author{E.~Lindenstrauss}
\address{E.L.: The Einstein Institute of Mathematics, Edmond J.\ Safra Campus, 
Givat Ram, The Hebrew University of Jerusalem, Jerusalem, 91904, Israel}
\email{elon@math.huji.ac.il}
\thanks{E.L.\ acknowledges support by ERC 2020 grant HomDyn (grant no.~833423).}

\author{A.~Mohammadi}
\address{A.M.: Department of Mathematics, University of California, San Diego, CA 92093}
\email{ammohammadi@ucsd.edu}
\thanks{A.M.\ acknowledges support by the NSF grants DMS-1764246 and 2055122.}

\author{Z.~Wang}
\address{Z.W.: Pennsylvania State University,
Department of Mathematics, University Park, PA 16802}
\email{zhirenw@psu.edu}
\thanks{Z.W.\ acknowledges support by the NSF grant  DMS-1753042.}

\author{L.~Yang}
\address{L.Y.: College of Mathematics, Sichuan University, Chengdu, Sichuan, 610000, China}
\email{lyang861028@gmail.com}
\thanks{L.Y.\ acknowledges support by the NSFC grant 12171338.}

\date{}

\begin{document}

\begin{abstract}
 We prove an effective version of the Oppenheim conjecture with a polynomial error rate.
The proof is based on an effective equidistribution theorem which in turn relies on recent progress towards restricted projection problem.   
\end{abstract}

\maketitle

\section{Introduction}

We prove the following theorem.

\begin{thm}\label{thm: Oppenhiem}
 Let $Q$ be an indefinite ternary quadratic form with $\det Q=1$. For all $R$ large enough, depending on $\norm{Q}$, and all $T\geq R^{\ref{a:Oppenheim}}$ at least one of the following holds.
\begin{enumerate}
    \item For every $s\in [-R^{\ref{k:Oppenheim}}, R^{\ref{k:Oppenheim}}]$, there exists a primitive vector $v\in \Z^3$ with $0<\norm{v}< T$ so that 
    \[
    \absolute{Q(v)-s}\leq R^{-\ref{k:Oppenheim}}
    \]
    \item There exists an integral quadratic form $Q'$ with $\norm{Q'}\leq R$ so that 
    \[
    \norm{Q-\lambda Q'}\leq  R^{\ref{a:Oppenheim}}(\log T)^{\ref{a:Oppenheim}} T^{-1}\qquad\text{where $\lambda=(\det Q')^{-1/3}$.}
    \]
\end{enumerate}
The constants $\consta\label{a:Oppenheim}$ and $\constk\label{k:Oppenheim}$ 
are absolute.
\end{thm}

Since algebraic numbers cannot be well approximated by rationals, one concludes the following corollary from Theorem~\ref{thm: Oppenhiem}. 

\begin{coro}\label{cor: Oppenheim}
Let $Q$ be a reduced, indefinite, ternary quadratic form which is not proportional to an integral form but has algebraic coefficients. Then for all $R$ large enough, depending on the degrees and heights of the coefficients of $Q$, and all $T\geq R^{\ref{a:Oppenheim 2}}$ we have the following: for any 
\[
s\in [-R^{\ref{k:Oppenheim 2}}, R^{\ref{k:Oppenheim 2}}]
\]
there exists a primitive vector $v\in \Z^3$ with $0<\norm{v}< T$ so that 
    \[
    \absolute{Q(v)-s}\leq R^{-\ref{k:Oppenheim 2}}.
    \]
The constants $\consta\label{a:Oppenheim 2}$ and $\constk\label{k:Oppenheim 2}$ 
are absolute.
\end{coro}

Theorem~\ref{thm: Oppenhiem} (and Corollary~\ref{cor: Oppenheim}) is an effective version of a celebrated theorem of Margulis~\cite{Margulis-Oppenheim}, see also~\cite{DM-Elemntary}. An effective version with a logarithmic rate was proved by the first named author and Margulis~\cite{LM-Oppenheim}. The proofs in~\cite{Margulis-Oppenheim, DM-Elemntary} and~\cite{LM-Oppenheim} rely on establishing a special case of Raghunathan's conjecture for unipotent flows --- albeit an effective version in the case of~\cite{LM-Oppenheim}. Raghunathan's conjecture in its full generality was proved by Ratner~\cite{Ratner-Acta, Ratner-measure, Ratner-topological}.  

Our proof of Theorem~\ref{thm: Oppenhiem} also is based on establishing an effective equidistribution theorem, as explicated below. Quantitative and effective versions of the aforementioned rigidity theorems in homogeneous dynamics have been sought after for some time, we refer to~\cite{LM-PolyDensity, LMW22, Yang-SL3} for a more detailed discussion of this problem and some recent progress.

Let 
\[
\sqf(\cox_1,\cox_2,\cox_3)=-2\cox_1\cox_3+\cox_2^2,
\]
and let $H=\SO(\sqf)^\circ\subset \SL_3(\R)=G$; 
note that $H\simeq\PSL_2(\R)$. 

For all $t,r\in\R$, let 
\[
a_t=\begin{pmatrix}
    e^t & 0 & 0\\
    0 & 1 & 0\\
    0 & 0 & e^{-t}
\end{pmatrix}
\quad\text{and}\quad u_r=\begin{pmatrix}
    1 & r & \frac{r^2}{2}\\
    0 & 1 & r\\
    0 & 0 & 1
\end{pmatrix}.
\]

Let $\Gamma\subset $ be a lattice. By Margulis' arithmeticity theorem, $\Gamma$ is arithmetic. $X=\SL_3(\R)/\Gamma$, and let $m_X$ denote the probability Haar measure on $X$. Let $\dist$ be the right invariant metric on $\SL_3(\R)$ which is defined using the Killing form and $\SO(3)$. This metric induces a metric $\dist_X$ on $X$, 
and natural volume forms on $X$ and its submanifolds.

\begin{thm}\label{thm: equidistribution 3}\label{thm: equidistribution n} 
For every $x_0\in X$ and large enough $R$ (depending explicitly on $x_0$), for any $T \geq R^{\ref{a:main-3-1}}$, at least one of the following holds.
\begin{enumerate}
\item For every $\varphi\in C_c^\infty(X)$, we have 
\[
\biggl|\int_0^1 \varphi(a_{\log T}u_rx_0)\diff\!r-\int \varphi\diff\!m_X\biggr|\leq \Sob(\varphi)R^{-\ref{k:main-3-1}}
\]
where $\Sob(\varphi)$ is a certain Sobolev norm. 
\item There exists $x\in X$ such that $Hx$ is periodic with $\vol(Hx)\leq R$, and 
\[
\dist_X(x,x_0)\leq R^{\ref{a:main-3-1}}(\log T)^{\ref{a:main-3-1}}T^{-1}.
\] 
\end{enumerate} 
The constants $\consta\label{a:main-3-1}$ and $\constk\label{k:main-3-1}$ are positive, and depend on $X$ but not on $x_0$. 
\end{thm}

The proof of Theorem~\ref{thm: equidistribution 3}, which is outlined in \S\ref{sec: the proof}, follows the same lines as the proof of~\cite[Thm.\ 1.1]{LMW22} where we replace the use of the projection theorem there with the recent work of Gan, Guo, and Wang~\cite{GGW}. 

\medskip

Using Theorem~\ref{thm: equidistribution 3} and arguments developed in~\cite{DM-Linearization, EMM-Upp}, we obtain the following strengthening of Theorem~\ref{thm: Oppenhiem}.

\begin{thm}\label{thm: lower bd 3}
Let $Q$ be an indefinite ternary quadratic form with $\det Q=1$. For all $R$ large enough, depending on $\norm{Q}$, and all $T\geq R^{\ref{a:Oppenheim L}}$ at least one of the following holds.
\begin{enumerate}
    \item Let $a<b$, then we have 
    \begin{multline*}
    \# \Bigl\{v\in\Z^3: \norm{v}\leq T, a\leq Q(v)\leq b\Bigr\}\geq \\  
    \mathsf C_Q(b-a)T + (1+\absolute{a}+\absolute{b})^NTR^{-\ref{k:Oppenheim L}}.
    \end{multline*}
    \item There exists an integral quadratic form $Q'$ with $\norm{Q'}\leq R$ so that 
    \[
    \norm{Q-\lambda Q'}\leq R^{\ref{a:Oppenheim L}}(\log T)^{\ref{a:Oppenheim L}} T^{-1}\qquad\text{where $\lambda=(\det Q')^{-1/3}$.}
    \]
\end{enumerate}
The constants $N$, $\consta\label{a:Oppenheim L}$, and $\constk\label{k:Oppenheim L}$ 
are absolute, and 
\[
\mathsf C_Q=\int_{L\cap B({\bf 0},1)}\frac{\diff\!\sigma}{\norm{\nabla Q}}
\]
where $B({\bf 0},1)$ is the ball of radius one in $\R^3$, $L=\{v: Q(v)=0\}$, and $\diff\!\sigma$ is the area element on $L$.  
\end{thm}

Note that the main term in part~(1) captures the asymptotic behavior of the volume of the solid given by $Q(v)=a$, $Q(v)=b$, and $\norm{v}\leq T$. We also remark that there is a dense family of irrational quadratic forms where the left side in part~(1) is $\gg T(\log T)^{1-\vare}$, see~\cite[\S 3.7]{EMM-Upp}, however, the forms in~\cite[\S 3.7]{EMM-Upp} are very well approximated by rational forms.

\section{The main steps in the proof}\label{sec: the proof}
The deduction of Theorem~\ref{thm: Oppenhiem} from Theorem~\ref{thm: equidistribution 3} is standard; and as it was mentioned, the deduction of Theorem~\ref{thm: lower bd 3} from Theorem~\ref{thm: equidistribution 3} follows from the arguments developed in~\cite{DM-Linearization, EMM-Upp},  see also~\cite{LMW23}. Thus, we turn to the proof of Theorem~\ref{thm: equidistribution 3}. 

\subsection*{The projection theorem of~\cite{GGW}}
Let $\gfrak=\Lie(G)$ and $\hfrak=\Lie(H)$, then  
$\gfrak=\hfrak\oplus\rfrak$
where $\rfrak$ is an irreducible representation of $\Ad(H)\simeq\PSL_2(\R)$ and $\dim\rfrak=5$. 

For every $s,z\in\R$ let 
\[
v_{s,z}=\begin{pmatrix}
    1 & -s & z\\
    0 & 1 & s\\
    0 & 0 & 1
\end{pmatrix}.
\]
Let $V=\{v_{s,z}: s, z\in\R\}$. Then $\Lie(V)=\rfrak^+$, where $\rfrak^+$ denotes the span of vectors in $\rfrak$ with positive $a_t$-weight. Moreover, if $W$ denotes the group of upper triangular unipotent matrices in $G$, then 
\[
W=UV.
\]

\begin{thm}[\cite{GGW}, Thm.~2.1]
\label{thm: projection}
Let $0<\alpha\leq 2$, and let  $0<\rhsc_1\leq1$. 
Let $\Theta\subset B_\rfrak(0,1)$ be a finite set satisfying the following 
\[
\#(B_\rfrak(w, \rhsc)\cap \Theta)\leq \egbd\cdot \rhsc^\alpha\cdot (\#\Theta)\qquad\text{for all $w$ and all $\rhsc\geq \rhsc_1$}
\]
where $\egbd\geq 1$.

Let $0<\pvare<10^{-4}\alpha$.
For every $\rhsc\geq \rhsc_1$, there exists a subset 
$J_{\rhsc}\subset [0,1]$ with $|[0,1]\setminus J_{\rhsc}|\ll \rhsc^{\pvare}$ so that the following holds. 
Let $r\in J_{\rhsc}$, then there exists a subset $\Theta_{\rhsc,r}\subset \Theta$ with 
\[
\#(\Theta\setminus \Theta_{\rhsc,r})\ll \rhsc^{\pvare}\cdot (\#\Theta)
\]
such that for all $w\in\Theta_{\rhsc, r}$, we have 
\[
\#\Bigl(\{w'\in\Theta: |\xi_{r}(w')-\xi_{r}(w)|\leq \rhsc\}\Bigr)\leq C\egbd\cdot \rhsc^{\alpha-O(\pvare)}\cdot (\#\Theta)
\] 
where $C\ll\pvare^{-\star}$, the implied constants are absolute, 
\[
\xi_{r}(w)= \Bigl(\Ad(r)w\Bigr)^+,
\]
and for every $w\in\rfrak$, $w^+$ denotes the orthogonal projection of $w$ to $\rfrak^+$. 
\end{thm}

As mentioned above $\rfrak$ is an irreducible $5$ dimensional representation of $\PSL_2(\R)$. Thus if we identify $\rfrak$ with $\R^5=\{(w_0,\ldots, w_4): w_i\in\R\}$ and the adjoint representation of $H$ on $\rfrak$ with the standard irreducible representation of $\SL_2(\R)$ on $\R^5$, then 
\[
\xi_{r}(w)=\Bigl(w_0+w_1r+w_2\tfrac{r^2}{2}+w_3\tfrac{r^3}{3!}+w_4\tfrac{r^4}{4!}, w_1+w_2r+w_3\tfrac{r^2}{2}+w_4\tfrac{r^3}{3!}\Bigr). 
\]
Theorem~\ref{thm: projection} thus follows from~\cite[Thm.~2.1]{GGW} applied with 
$n=5$, the curve $(r, r^2, r^3, r^4, r^5)$, and $m=2$.

\subsection*{Initial dimension}
As it was mentioned, the proof follows the same general steps as the proof of~\cite[Thm 1.1]{LMW22}. Indeed, assuming part~(2) in Theorem~\ref{thm: equidistribution n} does not hold, we first use an argument relying on Margulis functions for periodic orbits, to show that for $t_1=\log T-O(\log R)$ and all but a set with measure $\ll R^{-\star}$ of $r\in [0,1]$,  $x_1=a_{t_1}u_r x_0$ satisfies 
\be\label{eq: aways from periodic}
\dist_X(x,x_1)\gg R^{-D_0}
\ee
for all $x\in X$ with $\vol(Hx)\leq R$, where $D_0$ is absolute and the implied constants depend on $X$, see~\cite[Prop. 4.6]{LMW22}. 

\medskip

The rest of the argument is devoted to showing the following: let $x_1$ satisfy~\eqref{eq: aways from periodic}, then $\{a_{t-t_1}u_rx_1: u_r\in[0,1]\}$ will be $R^{-\star}$-equidistributed. 

To this end, we first use arithmeticity of $\Gamma$ to show that~\eqref{eq: aways from periodic} implies that points in $\{a_{\star \log R}u_rx_1:r\in[0,1]\}$ 
(possibly after removing an exceptional set of measure $\ll R^{-\star}$) are {\em separated transversal to} $H$. 

This follows from a closing lemma whose statement and the proof are similar to~\cite[Prop.\ 4.8]{LMW22}, mutatis mutandis; recall that $H\simeq \PSL_2(\R)$.

\subsection*{Improving the dimension}
We interpret separation transversal to $H$ which was obtained in the previous step as (small) positive dimension transversal to $H$ at controlled scales. In the next step, we use a localized Margulis function argument to improve this initial localized dimension transversal to $H$ to arrive at dimensions close to $2=\dim V$, recall that $W=UV$. 

This step is carried out as in~\cite[\S9--\S12]{LMW22}, see in particular~\cite[Prop.\ 10.1]{LMW22}. The main technical lemma in the argument is Lemma 9.1 in~\cite{LMW22} whose proof relies on tools from incidence geometry \'a la Wolff and Schlag. Our analysis here replaces the use of this incidence geometry technology with Theorem~\ref{thm: projection}. Let us now discuss this step in more details. 

Let $s=\star\log R$ and fix $0<\vare\leq 10^{-8}$. 
Let $\ell=0.01\vare s$, and let $\rwm$ be the probability measure on $H$ defined by 
\[
\rwm(\varphi)=\ave\varphi(a_{\ell}\uvk)\uvkd\qquad\text{for all $\varphi\in C_c(H)$.}
\] 

Our goal is to show that transversal dimension of $\{a_{s}u_rx_1: r\in[0,1]\}$, which was obtained in the initial dimension phase, 
will improve under convolution with $\nu$. It will be more convenient to {\em approximate} translations 
$\{a_su_rx_1: r\in[0,1]\}$ with sets which are a disjoint union of local $H$-orbits as we now define.  
Let $\kappa>0$ be much smaller than $\vare$, and put 
\[
\coneH=\Bigl\{u_{r'}^-: |r'|\leq e^{-\kappa s}\Bigr\}\cdot \Bigl\{a_\tau: |\tau|\leq e^{-\kappa s}\Bigr\}\Bigl\{u_r: |r|\leq e^{-\kappa s/10}\Bigr\}.
\]
Let $F\subset B_\rfrak(0, e^{-\kappa s})$ be a finite set with $\#F\geq \nuni^{s/2}$; 
let $y\in X$ with $\inj(y)>e^{-10\kappa s}$ (where $\inj$ denotes the injectivity radius of $y$), and  
put 
\be\label{eq: define cone intro}
\cone=\bigcup\coneH.\{\exp(w)y: w\in F\}.
\ee
We always assume $\cone$ is equipped with a probability measure $\mu_\cone$ 
which should be thought of as a weighted version of the uniform measure on $\cone$.

The conditional measures of $\mu_{\cone}$ along $\exp(B_\rfrak(0,1))$ are the objects of interest;  
the goal is to show that for $0<\alpha<2$, truncated $\alpha$-dimensional energies of these conditionals 
improve under convolution with $\nu$. This is achieved using the following localized and modified Margulis function. 

Let $0<\scmf\leq1/10$ (in our applications $\scmf=e^{-\sqrt\vare s}$, which is much smaller that $e^{-\ell}$). 
For every $(h,z)\in H\times \cone$, define 
\[
\margI_{\cone,\scmf}(h,z):=\Bigl\{w\in \rfrak: \|w\|<\scmf\,\inj(hz),\, \exp(w) h z\in h\cone.x\Bigr\};
\]
since $\mathsf E$ is bounded, $\margI_{\cone,\scmf}(h,z)$ is a finite set for all $(h,z)\in H\times\cone$.

For every $\trct\geq 1$, define  
$
\mfm_{\cone,\scmf,\trct}: H\times \cone\to [1,\infty)
$
as follows: 
\[
\text{if} \quad\#\margI_{\cone,\scmf}(h,z)\leq \trct,\quad \text{put}\quad \mfm_{\cone,\scmf,\trct}(h,z)=(\scmf\,\inj(hz))^{-\alpha}
\]
and if $\#\margI_{\cone,\scmf}(h,z)>\trct$, put
\[
\mfm_{\cone,\scmf,\trct}(h,z)=\min\left\{\sum_{w\in I}\|w\|^{-\alpha}: \begin{array}{c}I\subset I_{\cone,\scmf}(h,z)\text{ and }\\ \#(I_{\cone,\scmf}(h,z)\setminus I)=\trct\end{array}\right\}.
\]

The initial dimension phase produces sets $\cone$, see~\eqref{eq: define cone intro}, satisfying  
\[
\mfm_{\cone,\scmf,0}(e,z)\leq R^{D}
\] 
for some $D$ (depending only on $X$). Our goal is to show that
$\nu*\mu_\cone$ can be approximated (up to a set with measure $\ll R^{-\star}$) with sets $\cone_{\rm new}$ 
defined as in~\eqref{eq: define cone intro}, so that 
\be\label{eq: marg fun bound for new cones}
\mfm_{\cone_{\rm new},\scmf,e^{\ell/100}}(e,z)\leq e^{-\ell /2}R^{D}
\ee
This is achieved using Theorem~\ref{thm: projection} and the following observation: 
let $r\in[0,1]$ and let $w\in F$, then for any $F'\subset F$,   
\be\label{eq: expansion in plus}
\sum_{w'\in F'} \|\Ad(a_\ell u_r) (w-w')\|^{-\alpha}\ll \sum_{w'\in F'} e^{-\alpha\ell} \|(\Ad(u_r)(w-w'))^+\|^{-\alpha}.
\ee
Now using Theorem~\ref{thm: projection}, 
for all $r\in[0,1]$, but an exceptional set with measure $\ll R^{-\star}$, 
there is a set $F_r\subset F$ with $\#(F\setminus F_r)\ll R^{-\star}\cdot(\# F)$ so that for every $w\in F_r$ there is $F'\subset F_r$ satisfying the following 
\begin{align*}
  &\ \ \#(F_r\setminus F')\leq e^{\ell/100}\quad\text{and}\\
  &\sum_{w'\in F'} \|(\Ad u_r(w-w'))^+\|^{-\alpha}\ll e^{\ell/100} R^D.   
\end{align*}
Combining this and~\eqref{eq: expansion in plus} one proves~\eqref{eq: marg fun bound for new cones}, see~\cite[Lemma 9.1]{LMW22}. 

\medskip

That we can arrive at dimensions close to two relies on the fact that there are no intermediate obstructions and is made possible thanks to Theorem~\ref{thm: projection}. This appears to be analogous to the analysis in~\cite{LMW22, LM-PolyDensity}, where we increase the {\em invariance} group to the obtain the full horospherical subgroup $W$ in $G$, or at least an almost full dimensional subset of it. 
Note however that the codimension of $U$ in those cases was indeed one, thus no intermediate subgroup of $W$ existed.  
The fact that we avoid intermediate subgroups here is a new feature, when compared with the aforementioned works of Margulis, Dani, Ratner and others. E.g., in the work of Dani and Margulis~\cite{Margulis-Oppenheim, DM-Elemntary} one aims at increasing the invariance inductively, first obtaining the one dimensional group $Z(W)$ and then with the aid of the group $UZ(W)$ obtaining the full horospherical group $W$. We however obtain nearly full dimension in $W$ without passing through this intermediate step.

\subsection*{From large dimension to equidistribution}
The final step is to conclude equidistribution from this {\em high dimension}, as it is done in~\cite[Prop.\ 13.1]{LMW22}. As in loc.\ cit., this step is proved using an argument due to Venkatesh~\cite{Venkatesh-Sparse}, see~\cite[Prop.\ 5.2]{LMW22} which relies on the spectral gap of the ambient space, and a modification of the results concerning restricted projections. 

It is worth mentioning that unlike~\cite{LMW22} the horospherical subgroup $W=UV$ considered here is not commutative. In particular, conjugation by $a_t$ has two weights: the center direction $Z(W)$ is multiplied by $e^{2t}$ while other coordinates are multiplied by $e^t$. This, however, does not cause serious issues when adapting the proof of~\cite[Prop.\ 5.2]{LMW22} to the case at hand.

\bibliographystyle{halpha}
\bibliography{papers}

\end{document}